# Computational imaging of small-amplitude biperiodic surfaces with negative index material

BY YULIANG WANG

1. Research Center for Mathematics, Beijing Normal University, Zhuhai 519087, China

2. Guangdong Provincial Key Laboratory IRADS, BNU-HKBU United International College, Zhuhai 519087, China

*Email:* `yuliangwang@bnu.edu.cn`

**Abstract**

This paper presents an innovative approach to computational acoustic imaging of biperiodic surfaces, exploiting the capabilities of an acoustic superlens to overcome the diffraction limit. We address the challenge of imaging physical entities in complex environments by considering the partial differential equations that govern the physics and solving the corresponding inverse problem. We focus on imaging infinite rough surfaces, specifically 2D diffraction gratings, and propose a method that leverages the transformed field expansion (TFE). We derive a reconstruction formula connecting the Fourier coefficients of the surface and the measured field, demonstrating the potential for unlimited resolution under ideal conditions. We also introduce an approximate discrepancy principle to determine the cut-off frequency for the truncated Fourier series expansion in surface profile reconstruction. Furthermore, we elucidate the resolution enhancement effect of the superlens by deriving the discrete Fourier transform of white Gaussian noise. Our numerical experiments confirm the effectiveness of the proposed method, demonstrating high subwavelength resolution even under slightly non-ideal conditions. This study extends the current understanding of superlens-based imaging and provides a robust framework for future research.

## 1 Introduction

Computational imaging is an advanced method that employs computational algorithms to extract information about physical entities from emitted or scattered waves. In contrast to direct imaging, which relies on lens systems, computational imaging offers a versatile approach. Widely applied in ultrasound imaging, computed tomography, and magnetic resonance imaging, it enables accurate quantitative reconstructions and access to information beyond the reach of direct imaging techniques. However, both computational and direct imaging face a common limitation—the Abbe diffraction limit [1]—which dictates that the maximum achievable resolution is approximately half the wavelength of the wave used in the imaging process.

The field of metamaterials has witnessed significant progress since Victor Veselago introduced the concept of negative-index material [2]. Negative permittivity and permeability characterize these materials, resulting in optical properties contrary to traditional materials. John Pendry expanded this concept by proposing a superlens made from negative-index materials [3]. This superlens, utilizing a material with a negative refractive index, was theorized to amplify the evanescent field, achieving unlimited resolution in theory. This concept extended to acoustic waves [4], where a negative index involves both negative density and negative bulk modulus. Theoretical propositions were supported by numerous demonstrations of negative index acoustic metamaterials [5, 6, 7, 8]. Subsequently, the idea of an acoustic superlens was proposed and validated through various methods [9, 10, 11, 12].

However, when dealing with physical entities such as obstacles or inhomogeneous mediums, as in many applications, direct imaging by a superlens becomes challenging due to complicated scattered waves. We address this challenge through computational imaging, specifically by considering the partial differential equations underlying the physics and solving the corresponding inverse problem. This paper focuses on imaging infinite rough surfaces defined by biperiodic surfaces, known as 2D diffraction gratings. These gratings find applications not only in optics but also in acoustics for manipulating sound waves [13, 14, 15].

The forward scattering problem involves determining the diffracted field given a diffraction grating and incident field. The inverse scattering problem, on the other hand, aims to reconstruct the surface profile of the grating from measured data of the diffracted field for specific incident fields. Theoretical questions regarding the inverse diffraction problems, such as uniqueness and stability, have been extensively explored [16, 17, 18, 19, 20, 21, 22]. Computational methods often involve iterative algorithms [23, 24, 25, 26, 27, 28] or direct methods based on indicator functions [29, 30, 31, 32, 33, 34].

For surfaces with small amplitude, as assumed in this paper, the inverse problem can be linearized and solved using a method based on the transformed field expansion (TFE). Originally designed for direct rough surface scattering problems [35], TFE has been extended to various contexts, providing quantitative and computationally efficient solutions [36, 37, 38, 39, 40, 41, 42, 43, 44, 45, 46]. The small-amplitude assumption allows proving the uniqueness of the inverse problem for a single incident field, along with establishing convergence and error estimates for the computational method [47, 48].

This paper presents an imaging scheme where an acoustic superlens is placed over a biperiodic surface with small amplitude. A downward-propagating plane wave is incident upon the structure, and the field is measured on the superlens's top surface. Utilizing the transformed field expansion, we derive a reconstruction formula based on a simple relation connecting the Fourier coefficients of the surface and those of the measured field. We illustrate that achieving unlimited resolution is possible when the superlens is ideal. Furthermore, our numerical experiments indicate that high subwavelength resolution can be attained even for slightly non-ideal parameters.

Extending previous work [49], where we initially proposed imaging periodic surfaces with a superlens in the framework of inverse scattering, this paper addresses biperiodic surfaces with additional considerations. First, we introduce an approximate discrepancy principle to determine the cut-off frequency for the truncated Fourier series expansion in surface profile reconstruction, demonstrating its effectiveness through numerical experiments. Second, we explicitly derive the discrete Fourier transform of white Gaussian noise, elucidating the resolution enhancement effect of the superlens. Third, when comparing the imaging method's performance, we utilize signal-to-noise ratio rather than absolute noise level as a fair reference in practical applications.

The remainder of the paper is organized as follows. In Section 2, we establish the mathematical model for the forward diffraction problem. In Section 3, we apply the transformed field expansion to obtain zero-th and first-order closed-form solutions for the forward problem. These solutions are used in Section 4 to derive a reconstruction formula for the inverse problem and devise an approximate discrepancy principle for the cut-off frequency. We conduct comprehensive numerical experiments in Section 5 to test various aspects of the proposed method using smooth and non-smooth surface profiles. In Section 6, we conclude the paper with directions for future research.

## 2 Forward scattering problem

Consider the Helmholtz equation:
$$\nabla \cdot \left(\frac{1}{\rho} \nabla u\right) + \frac{\omega^2}{\kappa} u = 0, \tag{1}$$
which serves as the foundational partial differential equation for various wave phenomena. Here, $\omega$ represents the angular frequency arising from the implicit time dependence $e^{-i\omega t}$. In the context of acoustic waves, $\rho$ denotes the mass density, and $\kappa$ corresponds to the bulk modulus. In the case of electromagnetic waves, the variables $\rho$ and $\kappa$ respectively correspond to the electric permittivity and magnetic permeability, or vice versa, contingent upon the specific formulation of the physical model. Without loss of generality, we consider acoustic waves in the rest of the paper.

Let $\Gamma_f : z = f(x, y), (x, y) \in \mathbb{R}^2$ be a biperiodic surface, where
$$f(x + n_1 \Lambda_1, y + n_2 \Lambda_2) = f(x, y), \quad \forall (n_1, n_2) \in \mathbb{Z}^2$$

for some positive constants $\Lambda_1, \Lambda_2$. In addition, we assume $f$ is a small deformation of the plane $z = 0$, i.e.

$$f(x,y) = \varepsilon g(x,y), \tag{2}$$

where $0 < \varepsilon \ll 1$ is a prescribed constant, referred to as the deformation parameter in this paper.

Suppose $\rho = \kappa = 1$ in the half space $z > f(x,y)$. Let the incident wave be given by

$$u^{\text{in}} = e^{-i\omega z}, \quad z > f(x,y), \tag{3}$$

a plane wave propagating in the $-\hat{z}$ direction.

Assume $\Gamma_f$ is impenetrable so that a boundary condition can be imposed. For sound soft boundary, we have the Dirichlet boundary condition

$$u = 0 \quad \text{on } \Gamma_f.$$

For sound hard boundary, we have the Neumann boundary condition $\partial_\nu u = 0$. A Robin boundary condition can be used in general. In this paper, we focus on the Dirichlet boundary condition. The extensions to other boundary conditions can be considered routine.

Next we place a slab $\Omega : a < z < b$ above $\Gamma_f$, and its constituting medium has constant mass density $\rho$ and constant bulk modulus $\kappa$. We define the boundaries $\Gamma_b : z = b, \Gamma_a : z = a$ and the domains $\Omega : a < z < b, \Omega_f : f(x,y) < z < a$. A conceptual diagram is shown in Figure 1.

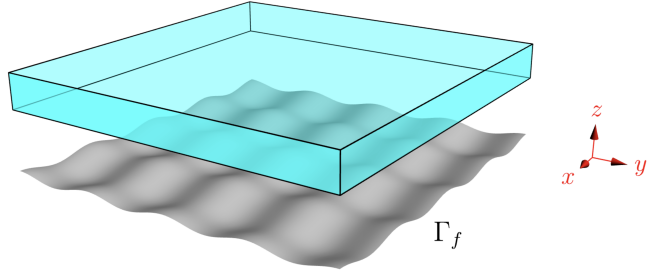

**Figure 1.** A conceptual diagram for the imaging problem.

By the continuity of the pressure field and the normal component of the velocity field across $\Gamma_a$ and $\Gamma_b$, we have the interface conditions

$$u^+ = u^-, \quad \frac{1}{\rho}\partial_z^+ u = \partial_z^- u \quad \text{on } \Gamma_a$$

$$u^+ = u^-, \quad \partial_z^+ u = \frac{1}{\rho}\partial_z^- u \quad \text{on } \Gamma_b \tag{4}$$

where the superscripts $+$ and $-$ denote the function values or partial derivatives taken from above and from below the interfaces, respectively.

Given the normal incident field (3), the total field $u$ is biperiodic in $x, y$ with the same periods of $f$. Let $u^{\text{sc}} = u - u^{\text{in}}$ denote the scattered field. Imposing the upward propagating radiation condition for $u^{\text{sc}}$ leads to the Rayleigh expansion

$$u^{\text{sc}}(x,y,z) = \sum_{\boldsymbol{n} \in \mathbb{Z}^2} u_{\boldsymbol{n}}^{\text{sc}}(b) e^{i[\boldsymbol{\alpha}_{\boldsymbol{n}} \cdot \boldsymbol{x} + \gamma_{\boldsymbol{n}}(z-b)]}, \quad z \geq b, \tag{5}$$

where $\boldsymbol{n} = (n_1, n_2) \in \mathbb{Z}^2, \boldsymbol{x} = (x,y) \in \mathbb{R}^2$,

$$\boldsymbol{\alpha}_{\boldsymbol{n}} = \left(\frac{2\pi n_1}{\Lambda_1}, \frac{2\pi n_2}{\Lambda_2}\right), \quad \gamma_{\boldsymbol{n}} = \sqrt{\omega^2 - |\boldsymbol{\alpha}_{\boldsymbol{n}}|^2}, \quad \text{Im}\,\gamma_{\boldsymbol{n}} \geq 0, \tag{6}$$

and $u_{\boldsymbol{n}}^{\text{sc}}(b)$ denotes the $\boldsymbol{n}$-th Fourier coefficients of $u^{\text{sc}}(\boldsymbol{x}, b)$, i.e.

$$u_{\boldsymbol{n}}^{\text{sc}}(b) = \frac{1}{\Lambda_1 \Lambda_2} \int_0^{\Lambda_2} \int_0^{\Lambda_1} u^{\text{sc}}(\boldsymbol{x}, b) e^{-i\boldsymbol{\alpha}_{\boldsymbol{n}} \cdot \boldsymbol{x}} dx dy.$$

Henceforward, we consistently employ the subscript $\boldsymbol{n}$ to signify the $\boldsymbol{n}$-th Fourier coefficient of a biperiodic function with respect to the variable $\boldsymbol{x}$. Note that the $\boldsymbol{n}$-th term in (5) is a propagating wave if $|\alpha_{\boldsymbol{n}}| < \omega$, and an evanescent wave if $|\alpha_{\boldsymbol{n}}| > \omega$. The evanescent wave is exponentially decreasing in amplitude, causing low signal-to-noise ratio when measured far from the surface. This can be understood as the physical origin of diffraction limit encountered in traditional imaging techniques.

Applying $\partial_z$ on both sides of (5) at $z = b$, we obtain

$$\partial_z^+ u^{\text{sc}} = T(u^{\text{sc}})^+ \quad \text{on } \Gamma_b, \tag{7}$$

where the operator $T$ is defined as

$$Tv = \sum_{\boldsymbol{n} \in \mathbb{Z}^2} i\gamma_{\boldsymbol{n}} v_{\boldsymbol{n}} e^{i\alpha_{\boldsymbol{n}} \cdot \boldsymbol{x}}$$

for any biperiodic function $v(\boldsymbol{x})$. Straightforward calculation shows

$$\partial_z^+ u^{\text{in}} = T(u^{\text{in}})^+ + \tau \quad \text{on } \Gamma_b, \tag{8}$$

where

$$\tau = -2i\omega e^{-i\omega b}. \tag{9}$$

Adding (7) to (8) and noting the continuity condition for $u$, we obtain

$$\partial_z^+ u = Tu + \tau \quad \text{on } \Gamma_b.$$

Combining this and the interface condition (4), we obtain the following boundary condition for the total field

$$\frac{1}{\rho}\partial_z^- u = Tu + \tau \quad \text{on } \Gamma_b.$$

We consolidate the derived differential equations, boundary conditions, and interface conditions to establish a boundary-interface value problem (BIVP) in $\Omega \cup \Omega_f$:

$$\begin{cases} \frac{1}{\rho}\partial_z u = Tu + \tau & \text{on } \Gamma_b, \\ \left[\Delta + \left(\frac{\rho}{\kappa}\right)\omega^2\right]u = 0 & \text{in } \Omega, \\ u^+ = u^-, \frac{1}{\rho}\partial_z^+ u = \partial_z^- u & \text{on } \Gamma_a, \\ (\Delta + \omega^2)u = 0 & \text{in } \Omega_f, \\ u = 0 & \text{on } \Gamma_f. \end{cases} \tag{10}$$

Given the surface profile $f$ and the incident field $u^{\text{in}}$, the forward scattering problem is to determine the total field $u$ for $z > f(x,y)$.

**Inverse scattering problem.** *To reconstruct the surface profile $f(x,y)$ given the incident field $u^{\text{in}}$ and the noisy measurement of total field $u$ on $\Gamma_b$.*

# 3 Transformed field expansion

The resolution of the inverse scattering problem inherently requires a thorough examination of the associated forward problem. As an initial step towards this objective, we derive a series solution for the BIVP as defined in (10), employing a transformed field expansion methodology.

## 3.1 Coordinate transformation

First, we make the following transformation of coordinates:

$$\tilde{x} = x, \tilde{y} = y, \tilde{z} = z, \quad (x,y,z) \in \bar{\Omega},$$

$$\tilde{x} = x, \tilde{y} = y, \tilde{z} = a\left[\frac{z - f(x,y)}{a - f(x,y)}\right], \quad (x,y,z) \in \bar{\Omega}_f,$$

This transformation maps the domain $\Omega_f$ to $\Omega_0 : 0 < a < b$, and the boundary $\Gamma_f$ to the flat surface $\Gamma_0 : z = 0$. The domain $\Omega$, the boundary $\Gamma_b$, and the interface $\Gamma_a$ remain fixed during this transformation.

We introduce the transformed function $\tilde{u}(\tilde{x}, \tilde{y}, \tilde{z}) = u(x, y, z)$, where $(x, y, z) \in \overline{\Omega \cup \Omega_f}$. After a meticulous yet straightforward calculation, we deduce from (10) that $\tilde{u}(\tilde{x}, \tilde{y}, \tilde{z})$, upon omitting the tildes from all variables, satisfies the equation

$$(c_1 \partial_{xx} + c_1 \partial_{yy} + c_2 \partial_{zz} - c_3 \partial_{xz} - c_4 \partial_{yz} - c_5 \partial_z + c_1 \omega^2) u = 0 \quad \text{in } \Omega_0, \tag{11}$$

where the coefficient functions $c_i$ are given by

$$\begin{aligned} c_1 &= (a - f)^2, \\ c_2 &= a^2 + (a - z)^2 |\nabla f|^2, \\ c_3 &= 2(a - z)(a - f) f_x, \\ c_4 &= 2(a - z)(a - f) f_y, \\ c_5 &= (a - z)[2|\nabla f|^2 + (a - f)\Delta f]. \end{aligned}$$

Furthermore, the interface condition in (10) is transformed to

$$u^+ = u^-, \quad \frac{1}{\rho}\left(1 - \frac{f}{a}\right)\partial_z^+ u = \partial_z^- u \quad \text{on } \Gamma_a. \tag{12}$$

The remaining equations and boundary conditions specified in (10) are invariant under the aforementioned transformation, retaining their original forms.

## 3.2 Asymptotic power series expansion

Under the small-amplitude assumption (2), we utilize either physical intuition or the transformed equation (11) to express the total field $u$ into the asymptotic power series expansion

$$u(x, y, z) = \sum_{k=0}^{\infty} \varepsilon^k u^{(k)}(x, y, z). \tag{13}$$

Substituting (13) into (11) and incorporating (2) into the expressions for $c_i$, followed by rearranging terms based on the powers of $\varepsilon$, yields the recursive system of equations

$$(\Delta + \omega^2) u^{(k)} = v^{(k)} \quad \text{in } \Omega_0,$$

where $v^{(0)} = 0$ and

$$v^{(1)} = a^{-1}[2g(\partial_{xx} + \partial_{yy}) + 2(a - z)(g_x \partial_{xz} + g_y \partial_{yz}) + (a - z)\Delta g \partial_z + 2\omega^2 g] u^{(0)}. \tag{14}$$

While the explicit forms for $v^{(k)}$ with $k \geqslant 2$ can be derived, they are deliberately omitted for the sake of brevity, as they are not relevant to the subsequent analysis.

Subsequently, by substituting (2) and (13) into the remaining equations in (10), we derive the following BIVP for $u^{(0)}$:

$$\begin{cases} \frac{1}{\rho}\partial_z u^{(0)} = T u^{(0)} + \tau & \text{on } \Gamma_b, \\ \left[\Delta + \left(\frac{\rho}{\kappa}\right)\omega^2\right] u^{(0)} = 0 & \text{in } \Omega, \\ [u^{(0)}]^+ = [u^{(0)}]^-, \frac{1}{\rho}\partial_z^+ u^{(0)} = \partial_z^- u^{(0)} & \text{on } \Gamma_a, \\ (\Delta + \omega^2) u^{(0)} = 0 & \text{in } \Omega_0, \\ u^{(0)} = 0 & \text{on } \Gamma_0, \end{cases} \tag{15}$$

and the following BIVP for $u^{(1)}$:

$$\begin{cases} \frac{1}{\rho}\partial_z u^{(1)} = T u^{(1)} & \text{on } \Gamma_b, \\ \left[\Delta + \left(\frac{\rho}{\kappa}\right)\omega^2\right] u^{(1)} = 0 & \text{in } \Omega, \\ [u^{(1)}]^+ = [u^{(1)}]^-, \frac{1}{\rho}\partial_z^+ u^{(1)} = \partial_z^- u^{(1)} + \frac{g}{\rho a}\partial_z^+ u^{(0)} & \text{on } \Gamma_a, \\ (\Delta + \omega^2) u^{(1)} = v^{(1)} & \text{in } \Omega_0, \\ u^{(1)} = 0 & \text{on } \Gamma_0. \end{cases} \quad (16)$$

## 3.3 Closed-form solutions

Since (15) consists of equations and boundary-interface conditions with constant coefficients in a rectangular domain, closed-form solutions become attainable. Applying the Fourier transform to (15) in $\boldsymbol{x}$ yields the following one-dimensional BIVP for $u_n^{(0)}(z)$:

$$\begin{cases} \frac{1}{\rho}\partial_z u_n^{(0)} = (i\gamma_n) u_n^{(0)} + \tau_n, & z = b, \\ (\partial_{zz} + \eta_n^2) u_n^{(0)} = 0, & a < z < b, \\ [u_n^{(0)}]^+ = [u_n^{(0)}]^-, \frac{1}{\rho}\partial_z^+ u_n^{(0)} = \partial_z^- u_n^{(0)}, & z = a, \\ (\partial_{zz} + \gamma_n^2) u_n^{(0)} = 0, & 0 < z < a, \\ u_n^{(0)} = 0, & z = 0, \end{cases} \quad (17)$$

where

$$\eta_n = \sqrt{\left(\frac{\rho}{\kappa}\right)\omega^2 - |\boldsymbol{\alpha}_n|^2}, \quad \text{Im}\,\eta_n \geqslant 0.$$

We retain the symbol $\partial$ for ordinary differentiation to maintain straightforward notations.

Clearly the general solution of (17) can be written as

$$u_n^{(0)}(z) = \begin{cases} A_n^{(0)} e^{i\eta_n z} + B_n^{(0)} e^{-i\eta_n z}, & a < z < b, \\ C_n^{(0)} e^{i\gamma_n z} + D_n^{(0)} e^{-i\gamma_n z}, & 0 < z < a, \end{cases} \quad (18)$$

with undetermined coefficients $A_n^{(0)}, B_n^{(0)}, C_n^{(0)}, D_n^{(0)}$. Applying the boundary and interface conditions in (17) yields the linear system of equations

$$M_n \begin{pmatrix} A_n^{(0)} \\ B_n^{(0)} \\ C_n^{(0)} \\ D_n^{(0)} \end{pmatrix} = \begin{pmatrix} \tau_n \\ 0 \\ 0 \\ 0 \end{pmatrix}, \quad (19)$$

where

$$M_n = \begin{pmatrix} i\psi_n e^{i\eta_n b} & -i\phi_n e^{-i\eta_n b} & 0 & 0 \\ e^{i\eta_n a} & e^{-i\eta_n a} & -e^{i\gamma_n a} & -e^{-i\gamma_n a} \\ \frac{i}{\rho}\eta_n e^{i\eta_n a} & -\frac{i}{\rho}\eta_n e^{-i\eta_n a} & -i\gamma_n e^{i\gamma_n a} & i\gamma_n e^{-i\gamma_n a} \\ 0 & 0 & 1 & 1 \end{pmatrix}, \quad (20)$$

and $\phi_n = \frac{\eta_n}{\rho} + \gamma_n$, $\psi_n = \frac{\eta_n}{\rho} - \gamma_n$ are introduced to simplify notations.

Direct calculation shows

$$\Sigma_n = \det(M_n) = e^{-i\gamma_n a}(e^{-i\eta_n h}\phi_n^2 - e^{i\eta_n h}\psi_n^2) + e^{i\gamma_n a}(e^{i\eta_n h} - e^{-i\eta_n h})\phi_n\psi_n, \tag{21}$$

where $h = b - a$ denotes the thickness of the superlens. Applying the Cramer's rule to (19) yields the solutions

$$\begin{cases} A_n^{(0)} = \dfrac{ie^{-i\eta_n a}\tau_n}{\Sigma_n}(\psi_n e^{-i\gamma_n a} - \phi_n e^{i\gamma_n a}), \\[4pt] B_n^{(0)} = \dfrac{ie^{i\eta_n a}\tau_n}{\Sigma_n}(\phi_n e^{-i\gamma_n a} - \psi_n e^{i\gamma_n a}), \\[4pt] C_n^{(0)} = -\dfrac{2i\eta_n \tau_n}{\rho \Sigma_n}, \\[4pt] D_n^{(0)} = -C_n^{(0)}. \end{cases} \tag{22}$$

Given that $\tau$ is a constant, we can immediately deduce that $\tau_0 = \tau$ and $\tau_n = 0$ for $n \neq 0$. Consequently, the coefficients $A_n^{(0)}, B_n^{(0)}, C_n^{(0)}$, and $D_n^{(0)}$ are all zero when $n \neq 0$. It follows from (18) that the zeroth order solution is given by

$$u^{(0)} = \begin{cases} A_0^{(0)} e^{i\eta_0 z} + B_0^{(0)} e^{-i\eta_0 z}, & a < z < b, \\ 2iC_0^{(0)} \sin(\gamma_0 z), & 0 < z < a. \end{cases} \tag{23}$$

Upon applying the Fourier transform to (16) and recognizing that $u^{(0)}$ is solely a function of $z$, we derive the one-dimensional BIVP for $u_n^{(1)}(z)$:

$$\begin{cases} \dfrac{1}{\rho}\partial_z u_n^{(1)} = (i\gamma_n) u_n^{(1)}, & z = b, \\[3pt] (\partial_{zz} + \eta_n^2) u_n^{(1)} = 0, & a < z < b, \\[3pt] [u_n^{(1)}]^+ = [u_n^{(1)}]^-, \dfrac{1}{\rho}\partial_z^+ u_n^{(1)} = \partial_z^- u_n^{(1)} + \dfrac{\partial_z^+ u^{(0)}}{\rho a} g_n, & z = a, \\[3pt] (\partial_{zz} + \gamma_n^2) u_n^{(1)} = v_n^{(1)}, & 0 < z < a, \\[3pt] u_n^{(1)} = 0, & z = 0, \end{cases} \tag{24}$$

By variation of parameters, we write the general solution of (24) as

$$u_n^{(1)}(z) = \begin{cases} A_n^{(1)} e^{i\eta_n z} + B_n^{(1)} e^{-i\eta_n z}, & a < z < b \\ C_n^{(1)} e^{i\gamma_n z} + D_n^{(1)} e^{-i\gamma_n z} + w_n(z), & 0 < z < a \end{cases} \tag{25}$$

where $A_n^{(1)}, B_n^{(1)}, C_n^{(1)}, D_n^{(1)}$ are undetermined coefficients, and

$$w_n(z) = \gamma_n^{-1} \int_0^z \sin[\gamma_n(z - z')] v_n^{(1)}(z') dz' \quad \gamma_n \neq 0. \tag{26}$$

Substituting (25) into the boundary and interface conditions in (24) leads to the following linear system

$$M_n \begin{pmatrix} A_n^{(1)} \\ B_n^{(1)} \\ C_n^{(1)} \\ D_n^{(1)} \end{pmatrix} = \begin{pmatrix} 0 \\ w_n(a) \\ w_n'(a) + \dfrac{\partial_z^+ u^{(0)}(a)}{\rho a} g_n \\ 0 \end{pmatrix},$$

where $M_n$ is given by (20). Solving the above system by the Cramer's rule, we obtain

$$A_n^{(1)} = \dfrac{2\phi_n e^{-i\eta_n b}}{\Sigma_n}\left[\gamma_n \cos(\gamma_n a) w_n(a) - \sin(\gamma_n a)\left(w_n'(a) + \dfrac{\partial_z^+ u^{(0)}(a)}{\rho a} g_n\right)\right]. \tag{27}$$

From (26), we have

$$w_n(a) = \frac{1}{\gamma_n}\int_0^a \sin[\gamma_n(a-z)]v_n^{(1)}(z)dz,$$
$$w_n'(a) = \int_0^a \cos[\gamma_n(a-z)]v_n^{(1)}(z)dz.$$

Upon substituting the above formulae into (27) and conducting a series of algebraic simplifications, we arrive at the following expression:

$$A_n^{(1)} = -\frac{2\phi_n e^{-i\eta_n b}}{\Sigma_n}\left[\int_0^a \sin(\gamma_n z)v_n^{(1)}(z)dz + \sin(\gamma_n a)\frac{\partial_z^+ u^{(0)}(a)}{\rho a}g_n\right]. \tag{28}$$

To determine $v_n^{(1)}$, we substitute (23) into (14) and apply the Fourier transform, yielding

$$v_n^{(1)}(z) = \frac{2i}{a}C_0^{(0)}\gamma_0[2\gamma_0\sin(\gamma_0 z) - |\alpha_n|^2(a-z)\cos(\gamma_0 z)]g_n, \quad 0 < z < a.$$

Upon direct calculation, we find

$$\int_0^a \sin(\gamma_n z)\sin(\gamma_0 z)dz = \frac{1}{2}\left\{\frac{\sin[(\gamma_n-\gamma_0)a]}{\gamma_n-\gamma_0} - \frac{\sin[(\gamma_n+\gamma_0)a]}{\gamma_n+\gamma_0}\right\},$$
$$\int_0^a \sin(\gamma_n z)(a-z)\cos(\gamma_0 z)dz = \frac{a\gamma_n}{\gamma_n^2-\gamma_0^2} - \frac{1}{2}\left\{\frac{\sin[(\gamma_n+\gamma_0)a]}{(\gamma_n+\gamma_0)^2} + \frac{\sin[(\gamma_n-\gamma_0)a]}{(\gamma_n-\gamma_0)^2}\right\}.$$

Combining these results and noting that $|\alpha_n|^2 = \gamma_0^2 - \gamma_n^2$, we deduce

$$\int_0^a \sin(\gamma_n z)v_n^{(1)}(z)dz = \frac{2i}{a}C_0^{(0)}\gamma_0[a\gamma_n - \sin(\gamma_n a)\cos(\gamma_0 a)]g_n.$$

Subsequently, from (23) and the interface condition in (15), we find

$$\partial_z^+ u^{(0)}(a) = \rho\partial_z^- u^{(0)}(a) = 2i\rho C_0^{(0)}\gamma_0\cos(\gamma_0 a).$$

Inserting the previous two expressions in (28), we obtain

$$A_n^{(1)} = -\frac{4i}{\Sigma_n}C_0^{(0)}\gamma_0\gamma_n\phi_n e^{-i\eta_n b}g_n. \tag{29}$$

Following analogous steps, we find

$$B_n^{(1)} = -\frac{4i}{\Sigma_n}C_0^{(0)}\gamma_0\gamma_n\psi_n e^{i\eta_n b}g_n. \tag{30}$$

The coefficients $C_n^{(1)}$ and $D_n^{(1)}$ can also be determined, but they not needed in the sequel and are therefore omitted.

Substituting (29) and (30) into (25), and evaluating the expression at $z = b$, we derive

$$u_n^{(1)}(b) = -\frac{8i}{\rho\Sigma_n}C_0^{(0)}\gamma_0\gamma_n\eta_n g_n.$$

Finally, incorporating the expression of $C_0^{(0)}$ from (22), we arrive at a pivotal identify of this paper:

$$g_n = s_n u_n^{(1)}(b), \tag{31}$$

where

$$s_n = -\frac{\rho^2 \Sigma_0 \Sigma_n}{16\tau\gamma_0\eta_0\gamma_n\eta_n} \tag{32}$$

serves as the scaling factors that links the Fourier coefficients of the surface profile $g(x)$ with the linear approximation of $u(x,b)$.

**Remark 1.** Should the conditions $\gamma_n = 0$ or $\eta_n = 0$ arise for some $n \in \mathbb{Z}^2$, indicating the occurrence of resonance, the solutions will manifest in distinct forms. While these special cases warrant separate treatment, they are omitted from the present study for the sake of brevity and focus.

## 4 Inverse scattering problem

Recalling the inverse scattering problem, the objective is to reconstruct $f^\delta$ from the perturbed measurement data $u^\delta(x,b) = u(x,b) + \delta(x)$, where $\delta(x)$ represents additive measurement noise. By substituting $u(x,b)$ with the asymptotic expansion (13), we arrive at

$$u^\delta(x,b) = u^{(0)}(x,b) + \varepsilon u^{(1)}(x,b) + r(x,b) + \delta(x), \tag{33}$$

where $r(x,b) = \sum_{k=2}^{\infty} \varepsilon^k u^{(k)}(x,y)$ denotes the remainder term in the expansion. Combining (33) with (31), we deduce

$$f_n = s_n [u_n^\delta(b) - u_n^{(0)}(b) - r_n(b) - \delta_n]. \tag{34}$$

It is noteworthy that $u^{(0)}$, explicitly provided by (23), is independent of $f$.

Upon dropping the error terms, we define

$$f_n^\delta = s_n [u_n^\delta(b) - u_n^{(0)}(b)], \tag{35}$$

which serves as an approximate representation of $f_n$. Finally, we reconstruct the surface profile function $f$ using the truncated Fourier series expansion

$$f^{\delta,N}(x) = \text{Re} \sum_{\|n\|_\infty \leqslant N} f_n^\delta e^{i\alpha_n \cdot x}, \tag{36}$$

where $\|n\|_\infty = \max(n_1, n_2)$ and $N$ signifies the cut-off frequency. It should be noted that alternative norms for $n$ may also be employed in this context.

Combining (34)–(36), we find that the total reconstruction error $f^{\delta,N} - f$ can be decomposed as

$$f^{\delta,N} - f = E_1 + E_2 + E_3,$$

where each component is defined as

$$E_1 = \text{Re} \sum_{\|n\|_\infty \leqslant N} s_n r_n(b) e^{i\alpha_n \cdot x}, \quad E_2 = \text{Re} \sum_{\|n\|_\infty \leqslant N} s_n \delta_n e^{i\alpha_n \cdot x}, \quad E_3 = \sum_{\|n\|_\infty > N} f_n e^{i\alpha_n \cdot x}.$$

It is worth noting the origins of each error term.

1. $E_1$ emerges from dropping the higher order terms in the asymptotic expansion (13), i.e. approximating the nonlinear problem with a linear one.

2. $E_2$ is associated with the measurement noise $\delta(x)$.

3. $E_3$ results from the truncation of the Fourier series representation of $f$, which depends on the cut-off frequency $N$ and the spectral properties of $f$.

Next, we introduce a discrepancy principle to determine the cut-off frequency $N$, which plays the role of a regularization parameter for the inverse problem. To this end, let $u^{\delta,N}$ denote the solution to the forward scattering problem corresponding to the surface profile $f^{\delta,N}$.

Neglecting higher-order terms in the asymptotic expansion for $u^{\delta,N}$ and employing the approximation $f_n^{\delta,N} \approx f_n^\delta$, we arrive at the expression

$$u^{\delta,N}(\boldsymbol{x},b) \approx u^{(0)}(\boldsymbol{x},b) + \sum_{\|\boldsymbol{n}\|_\infty \leqslant N} s_n^{-1} f_n^\delta e^{i\boldsymbol{\alpha}_n \cdot \boldsymbol{x}}$$
$$= u^{(0)}(\boldsymbol{x},b) + \sum_{\|\boldsymbol{n}\|_\infty \leqslant N} [u_n^\delta(b) - u_n^{(0)}(b)] e^{i\boldsymbol{\alpha}_n \cdot \boldsymbol{x}}.$$

It follows that

$$u^\delta(\boldsymbol{x},b) - u^{\delta,N}(\boldsymbol{x},b) \approx \sum_{\|\boldsymbol{n}\|_\infty > N} [u_n^\delta(b) - u_n^{(0)}(b)] e^{i\boldsymbol{\alpha}_n \cdot \boldsymbol{x}} =: R^{\delta,N},$$

where $R^{\delta,N}$ denotes an approximation of the residual.

We propose the following approximate discrepancy principle of Morozov for the cut-off frequency $N$.

**Discrepancy principle.** *For any given noise $\delta$, select $N$ as the smallest integer satisfying*

$$\|R^{\delta,N}\|_2 < c\|\delta\|_2, \tag{37}$$

*where $c > 0$ is a prescribed constant.*

It is noteworthy that $\|R^{\delta,N}\|_2$ is non-increasing with respect to $N$. Consequently, for a fixed constant $c$, a higher noise level $\|\delta\|_2$ necessitates a smaller value of $N$. Given that $R^{\delta,N}$ is an approximate, rather than exact, residual, it is possible that $\|R^{\delta,N}\|_2 < \|\delta\|_2$. This permits the utilization of $c \leqslant 1$, a departure from the classical discrepancy principle of Morozov [50]. We employ this approximate discrepancy principle in the numerical experiments conducted in Section 5, and the outcomes are both satisfactory and consistent.

Rigorous theoretical examination of the proposed discrepancy principle, as well as the formal establishment of the well-posedness for both the direct and inverse problems, and the convergence and error analysis of the numerical scheme under consideration, fall outside the focus of the present paper. Those topics, given their interconnected and complex nature, are best addressed collectively in a dedicated theoretical study, which we hope to report in a future paper.

Next, we investigate the numerical method in some special cases regarding the superlens.

## 4.1 Without superlens

Consider the case when the slab is absent, i.e. $\rho = \kappa = 1$. Then we have

$$\eta_n = \gamma_n, \quad \eta_0 = \gamma_0 = \omega, \quad \Sigma_n = 4\gamma_n^2 e^{-i\gamma_n b}, \quad \Sigma_0 = 4\omega^2 e^{-i\omega b}$$

Substituting the above results and (9) into (32), we derive a simplified expression for the scaling factor:

$$s_n = \frac{e^{-i\gamma_n b}}{2i\omega}, \tag{38}$$

From the definition of $\gamma_n$ in (6), we make the following distinctions.

i. If $|\alpha_n| \leqslant \omega$, then $\gamma_n \in \mathbb{R}$, corresponding to propagating waves in the Rayleigh expansion (5) for the scattered wave. In this case, $|s_n| = 1/(2\omega)$ is constant in terms of $\boldsymbol{n}$. It follows from (31) that the reconstruction of the corresponding frequency modes $f_n$ are equally stable for those values of $\boldsymbol{n}$.

ii. If $|\alpha_n| > \omega$, then $\gamma_n = i\sqrt{|\alpha_n|^2 - \omega^2}$, corresponding to evanescent waves in (5). In this case, $|s_n| = \frac{1}{2\omega} e^{b\sqrt{|\alpha_n|^2 - \omega^2}} \to \infty$ exponentially as $\|\boldsymbol{n}\| \to \infty$. Consequently, it is increasingly unstable to reconstruct the higher frequency modes $f_n$ when the measurement data is contaminated by noise with flat power spectral density across all frequencies (e.g. additive white Gaussian noise, see Section 6).

Note that the number of $\boldsymbol{n}$ satisfying $|\alpha_{\boldsymbol{n}}| \leqslant \omega$ is finite and roughly proportional to $\omega$. Thus the number of frequency modes of $f$ amenable to stable reconstruction is proportional to $\omega$, which elucidates the diffraction limit inherent in traditional far-field imaging methodologies.

## 4.2 With superlens

Next we consider the case where $\rho, \kappa$ are both negative. First we consider the special case when $\rho = \kappa = -1$, corresponding to the superlens proposed in [3]. Then we have

$$\eta_{\boldsymbol{n}} = \gamma_{\boldsymbol{n}}, \quad \eta_0 = \gamma_0 = \omega, \quad \Sigma_{\boldsymbol{n}} = -4\gamma_{\boldsymbol{n}}^2 e^{-i\gamma_n(a-h)}, \quad \Sigma_0 = -4\omega^2 e^{-i\omega(a-h)}.$$

Substituting these results into (32), we derive

$$s_{\boldsymbol{n}} = \frac{e^{2i\omega h}}{2i\omega} e^{-i\gamma_n(a-h)}. \tag{39}$$

Taking the absolute value of $s_{\boldsymbol{n}}$, we find $a-h$ in (39) plays the role of $b$ in (38). Thus we may call $a-h$ the "effective measuring distance" in this case. Physically speaking, $a-h$ is approximately the acoustic path length as the wave travels from $\Gamma_f$ to $\Gamma_b$, where $a$ is approximately the path length from $\Gamma_f$ to $\Gamma_a$, and $-h$ is the acoustic path length from $\Gamma_a$ to $\Gamma_b$. For evanescent waves, this implies they decay exponentially from $\Gamma_f$ to $\Gamma_a$, and subsequently grows exponentially from $\Gamma_a$ to $\Gamma_b$. In essence, this observation underscores the role of the superlens as an effective amplifier of evanescent waves.

Let us take $a = h$, so that (39) reduces to

$$s_{\boldsymbol{n}} = \frac{e^{2i\omega h}}{2i\omega}. \tag{40}$$

Thus $|s_{\boldsymbol{n}}| = (2i\omega)^{-1}$ is constant with respect to $\boldsymbol{n}$. Consequently, all the frequency modes of $f$ are equally amenable to stable reconstructions. In this context, we achieved unlimited resolution within the realm of linear approximation. Notably, this outcome serves to substantiate the findings originally presented in [3], albeit in the context of imaging periodic surfaces. Alternatively, we can advance the argument that when $a = h$, it signifies a scenario wherein the effective measuring distance reduces to zero. Additionally, the condition $a = h$ can be interpreted as inducing a zero effective measuring distance, which suggests that no information pertaining to the surface is lost in the reconstruction process.

In practical applications, achieving materials with precisely defined effective properties may pose significant challenges. Specifically, the effective density $\rho$ and bulk modulus $\kappa$ may manifest non-zero imaginary components, attributable to diffusive scattering losses within the composite metamaterial [4]. In our numerical experiments, we employ $\rho = -1 + i\mu_1$ and $\kappa = -1 + i\mu_2$, wherein $\mu_1$ and $\mu_2$ are small, positive constants. Under these conditions, the expression (40) remains approximately valid, thereby leading us to anticipate the persistence of resolution enhancement effects.

# 5 Numerical experiments

We conduct numerical experiments to investigate the effectiveness of the proposed reconstruction formula and discrepancy principle. For simplicity, we fix $\Lambda_1 = \Lambda_2 = 1$ throughout the experiments. Let $\lambda = 2\pi/\omega$ denote the wavelength in the free space and we fix $\lambda = 1.1$.

**Remark 2.** The proposed method presumes prior knowledge of the surface's periodicities $\Lambda_1, \Lambda_2$, an assumption that may not be applicable in practical contexts. To circumvent this limitation, one could formulate an inverse problem where $\Lambda_1$ and $\Lambda_2$ are treated as unknown variables to be ascertained. These parameters could be reconstructed either independently or in conjunction with the surface's profile.

## 5.1 Discretization and implementation

For a specified profile function $f$, we solve the forward scattering problem utilizing the Perfectly Matched Layer (PML) technique in conjunction with the Finite Element Method (FEM). In the implementation of the FEM solver, cubic Lagrange elements are employed. Moreover, the mesh generator is configured such that the maximum element size is restricted to $\lambda/32$ on the boundaries $\Gamma_a$ and $\Gamma_b$, and $\lambda/8$ for other regions within the mesh.

The Finite Element Method (FEM) solution for the profile function $f$ is subjected to interpolation to yield the total field $u(x_i, b)$ on a uniform rectangular grid $x_i$ in the domain $[0,1]^2$. Here, $i = (i_1, i_2)$ and $x_i = \left(\frac{i_1}{I_1}, \frac{i_2}{I_2}\right)$ with $0 \leq i_1 < I_1$ and $0 \leq i_2 < I_2$, where $I_1$ and $I_2$ are positive integers.

Subsequently, the field value undergoes a perturbation by a noise component to synthesize measurement data. The characteristics of the noise are application-dependent; in the realm of imaging sciences, additive white Gaussian noise is frequently employed. The perturbed data is thus expressed as

$$u^\delta(x_i, b) = u(x_i, b) + \delta_i,$$

where both the real and imaginary components of $\delta_i$ are independent of $u(x_i, b)$ and adhere to an independent and identically distributed Gaussian distribution with zero mean and a standard deviation $\sigma$. To ensure reproducibility of the results, the seed for the random number generator is fixed.

The Signal-to-Noise Ratio (SNR) is defined in terms of the data energy, formulated as:

$$\text{SNR} = \frac{\|u(\cdot, b)\|_{l^2}}{\|\delta\|_{l^2}}.$$

Finally, it is noteworthy that $u_0^{(0)}(b) = A_0^{(0)} e^{i\eta_0 b} + B_0^{(0)} e^{-i\eta_0 b}$ and $u_n^{(0)}(b) = 0$ when $n \neq 0$. Utilizing these data, we compute $f_0^\delta$ from (35) and ultimately derive a discrete version of the reconstructed profile function $f^{\delta, N}$ in accordance with (36).

In general, the function $u(x, b)$ is not band limited, so its Fourier coefficients cannot be computed exactly from finitely many sample points, especially when it is perturbed by noise. However, if the field is sufficiently smooth, then its Fourier coefficients decays rapidly and the low-frequency modes $u_n^\delta(b)$ can be accurately approximated by the discrete Fourier transform (DFT)

$$U_n^\delta = \frac{1}{I_1 I_2} \sum_{i_1=0}^{I_1-1} \sum_{i_2=0}^{I_2-1} e^{-2\pi i \left(\frac{n_1 i_1}{I_1} + \frac{n_2 i_2}{I_2}\right)} u^\delta(x_i, b).$$

Without loss of generality, we take $I_1 = I_2 = I$ and refer to $I$ as the sample rate. For more accurate results, $I$ should be as large as possible. Moreover, it follows from the Appendix that the DFT of $\delta(x_i)$ follows i.i.d Gaussian distribution with deviation inversely proportional to $I$. Consequently, greater sample rate yields greater SNR for a given noise distribution. On the other hand, $I$ should be as small as possible for easier and faster measurement in practical scenario. As a guiding principle, the Petersen–Middleton theorem dictates that $I > 2N$ if accurate reconstructions are sought up to the cut-off frequency $N$. In the sequel, we fix $I = 99$ in all numerical experiments.

## 5.2 Surface profiles

We carry out detailed experimentation on three distinct surface profiles, each progressively more challenging than the last. The first profile function is given by $g(x) = p(x) + p(y)$, where

$$p(t) = \frac{1}{4}[0.5 + \sin(2\pi t) + \cos(4\pi t) + \sin(6\pi t)].$$

A pseudocolor plot of the profile function in $[0,1]^2$ is shown in Figure 2 (a). Note that this profile function is smooth and consists of nonzero Fourier modes $g_n$ with $\|n\|_\infty \leq 3$, which is seen from the pseudocolor plot of $|g_n|$ in Figure 2 (d).

The second profile is the periodic extension of $q(8x-4, 8y-4), (x,y) \in [0,1]^2$, where

$$q(s,t) = 0.3(1-s)^2 e^{-s^2-(t+1)^2} - (0.2s - s^3 - t^5)e^{-s^2-t^2} - 0.03 e^{-(s+1)^2-t^2}.$$

This profile function is better illustrated by a surface plot, which is shown in Figure 2 (b). This profile function is numerically smooth but apparently consists of infinitely many Fourier modes. So the spectral cut-off error never vanishes but decreases fast as $\|\boldsymbol{n}\|_\infty \to \infty$, which is seen from the pseudocolor plot of $|g_{\boldsymbol{n}}|$ in Figure 2 (e).

The third profile is non-smooth, piecewise constant function created from the gray-scale image as shown in Figure 2 (c). We set $g(x,y)$ to be the periodic extension of the indicator function for the image. The points of jump discontinuities are connected by lateral surfaces to create the final surface profile. Although the derivation of our method assumes $g$ to be twicely differentiable at least, we still expect the method to work even in this case in view of the Fourier approximation theory. On the other hand, we anticipate more challenges in achieving accurate reconstruction in view of the Gibbs phenomenon and the slower decay rate of the Fourier modes comparing with the previous two functions, which is seen from the pseudocolor plot of $|g_{\boldsymbol{n}}|$ in Figure 2 (f).

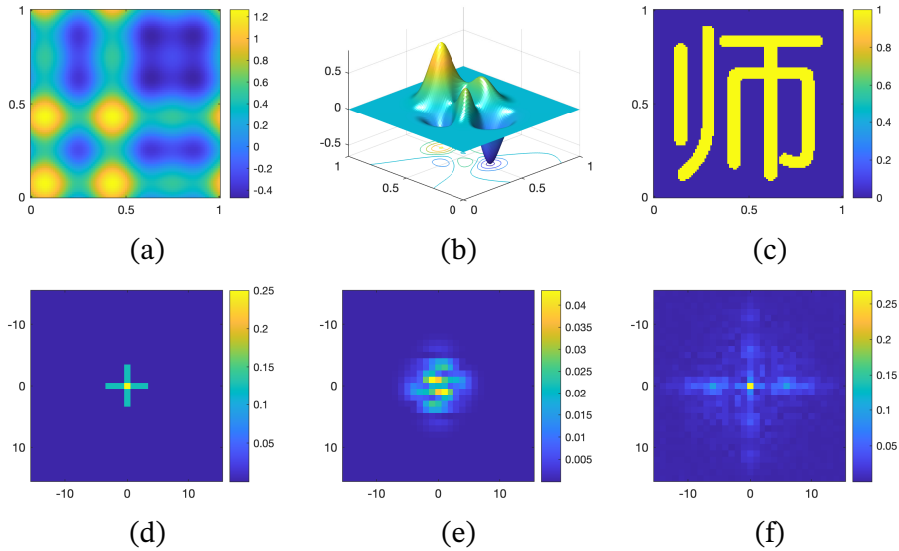

**Figure 2.** Ground truth of surface profiles and their Fourier coefficients. (a) pseudocolor plot of the first profile function; (b) surface plot of the second profile function; (c) pseudocolor plot of the third profile function. (d)–(f): pseudocolor plots of $|g_{\boldsymbol{n}}|$ for the profile functions $g$ in (a)–(c).

## 5.3 Experiment 1

We first consider the profile function shown in Figure 2 (a). With the parameter values given in Table 1 (row 1), we compute the reconstructions using cut-off frequencies $0 \leq N \leq 8$ and plot the results in Figure 3. We observe the Fourier modes of the surface profile is cumulatively reconstructed from low to high frequencies as $N$ increases. For this profile, $N = 3$ is sufficient but the reconstructions stays stable for several values greater than 3, demonstrating the robustness with respect to $N$. Certainly, increasing $N$ further would eventually leads to deteriorated reconstructions.

|  | $\varepsilon$ | $a$ | $h$ | $\rho$ | $\kappa$ | $\sigma$ | SNR |
|---|---|---|---|---|---|---|---|
| Figure 3 | 0.001 | 0.1 | 0.1 | $-1+0.01i$ | $-1+0.01i$ | 0.005 | 10.9 |
| Figure 4 | 0.001 | 0.1 | 0.1 | $-1+0.01i$ | $-1+0.01i$ | 0.01 | 5.5 |
| Figure 5 | 0.001 | 0.1 | 0.1 | $-1+0.01i$ | $-1+0.01i$ | 0.02 | 2.8 |

**Table 1.** Parameter values used in Experiment 1.

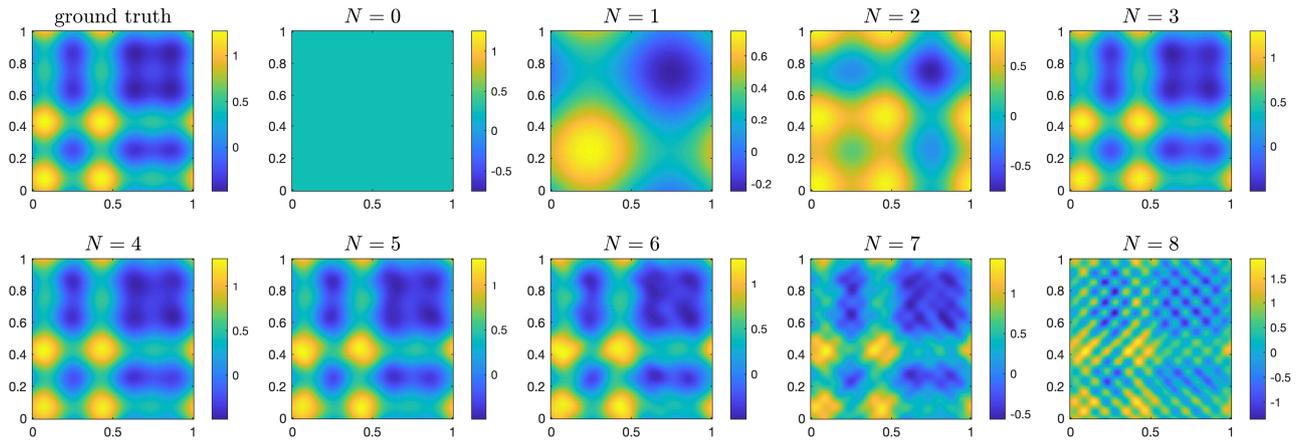

**Figure 3.** Pseudocolor plots of the reconstructions for the first surface profile with parameter values given in Table 1 (row 1) and different cut-off frequencies $N$.

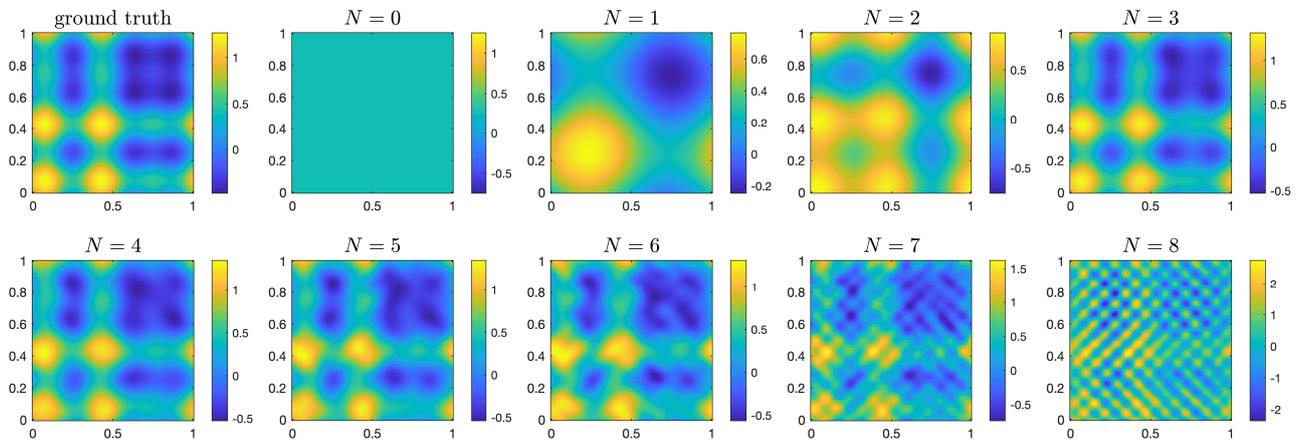

**Figure 4.** Pseudocolor plots of the reconstructions for the first surface profile with parameter values given in Table 1 (row 2) and different cut-off frequencies $N$.

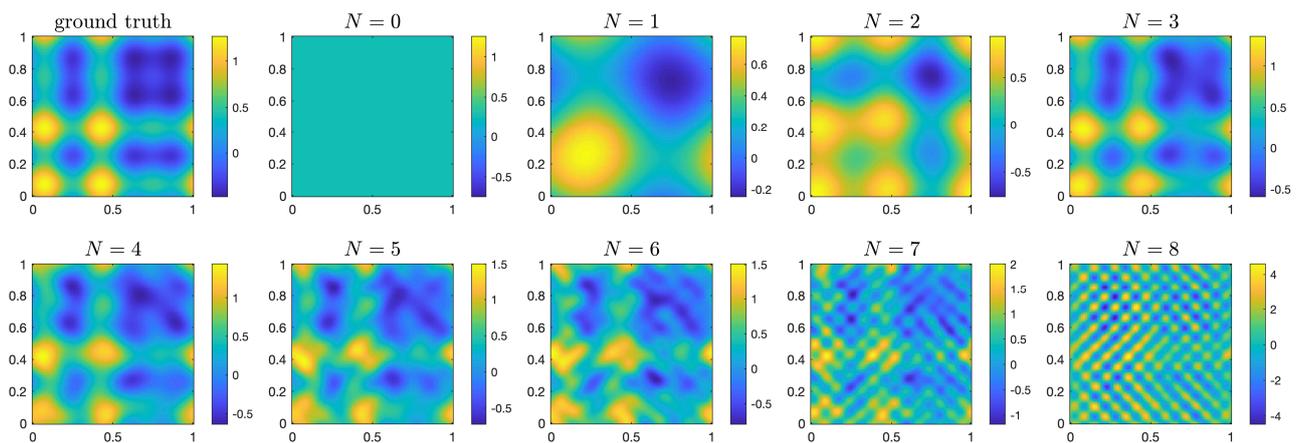

**Figure 5.** Pseudocolor plots of the reconstructions for the first surface profile with parameter values given in Table 1 (row 3) and different cut-off frequencies $N$.

Next we increase the noise level by setting $\sigma = 0.01$, leading to a signal-to-noise ratio SNR = 5.5. The results are shown in Figure 4. Observe the reconstructions exhibit minor difference with those in Figure 3 for small $N$ even with the much greater noise level, demonstrating high tolerance of the method with respect to measurement noise. On the other hand, the reconstructions deteriorate quicker than those in Figure 3 as $N$ increases further due to the greater noise level. As comparison, Figure 5 shows the results with $\sigma = 0.02$ and SNR = 2.8.

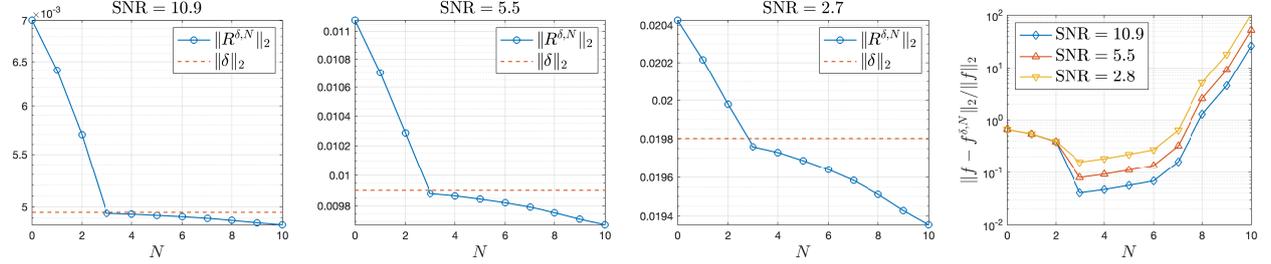

**Figure 6.** Row 1: residual $\|R^{\delta,N}\|_2$ (solid blue line) and absolute noise level $\|\delta\|_2$ (dashed red line). Row 2: relative error $\|f - f^{\delta,N}\|_2/\|f\|_2$. Column 1–3: with parameter values given in Table 1.

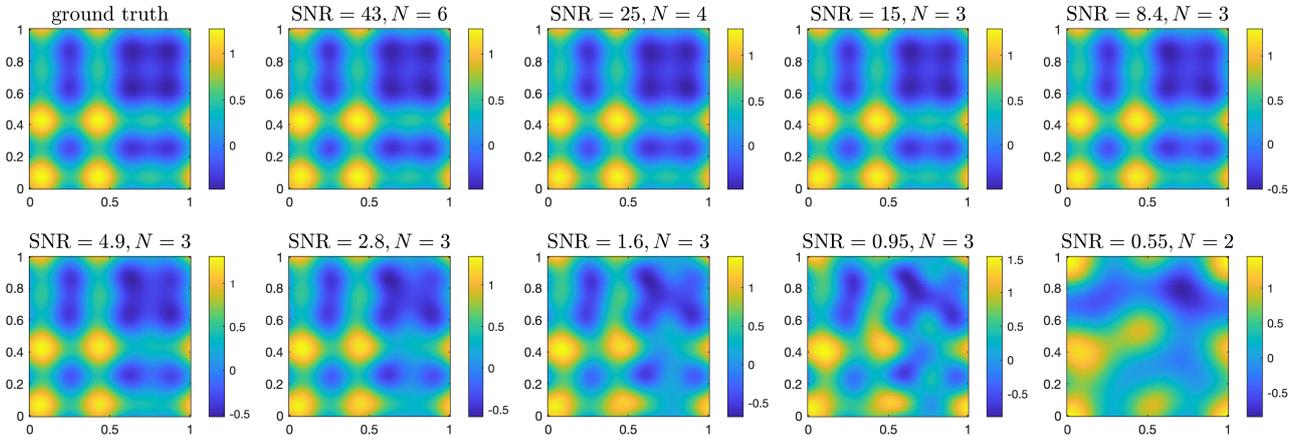

**Figure 7.** Pseudocolor plots of the reconstructions for the first surface profile using parameter values given in Table 1 (row 1), but with different values of SNR and the corresponding cut-off frequency $N$ determined from the discrepancy principle.

We may leverage the discrepancy principle (37) to determine a suitable cut-off frequency $N$. Letting $c = 1$, we plot the residual $\|R^{\delta,N}\|_2$ against $N$ and the absolute noise level $\|\delta\|_2$ for the previous experiments in Figure 6 (row 1). According to the discrepancy principle, we choose the smallest $N$ such that $\|R^{\delta,N}\|_2 \leqslant \|\delta\|_2$, leading to $N = 3$ for all three cases. This choice is optimal for this experiment as seen from the relative error $\|f - f^{\delta,N}\|_2/\|f\|_2$ of the reconstructions depicted in Figure 6 (row 2). In terms of the wavelength $\lambda = 1.1$, we obtain subwavelength resolutions even with very low signal-to-noise ratio using the proposed discrepancy principle.

Figure 7 presents reconstructions conducted with the parameter values specified in Table 1 (row 1), but with different values of SNR and the corresponding cut-off frequency $N$ determined from the discrepancy principle. Remarkably, the method exhibits robustness over a substantial SNR span, ranging from as low as 0.55 to as high as 44.

## 5.4 Experiment 2

We next consider the surface profile shown in Figure 2 (b). With the parameter values given in Ta-ble 2 (row 1), we compute the reconstructions using cut-off frequencies $0 \leqslant N \leqslant 8$ and plot the results in Figure 8. We again observe the Fourier modes of the surface profile is cumulatively reconstructed from low to high frequencies as $N$ increases at the beginning and it is increasingly unstable to reconstruct higher frequency modes.

|  | $\varepsilon$ | $a$ | $h$ | $\rho$ | $\kappa$ | $\sigma$ | SNR |
|---|---|---|---|---|---|---|---|
| Figure 8 | 0.001 | 0.1 | 0.1 | $-1+0.01i$ | $-1+0.01i$ | 0.005 | 9.4 |
| Figure 9 | 0.001 | 0.1 | 0.1 | $-1+0.001i$ | $-1+0.001i$ | 0.0009 | 9.3 |
| Figure 11 | 0.001 | 0.1 | 0.1 | 1 | 1 | 1 | 9.3 |

**Table 2.** Parameter values used in Experiment 2.

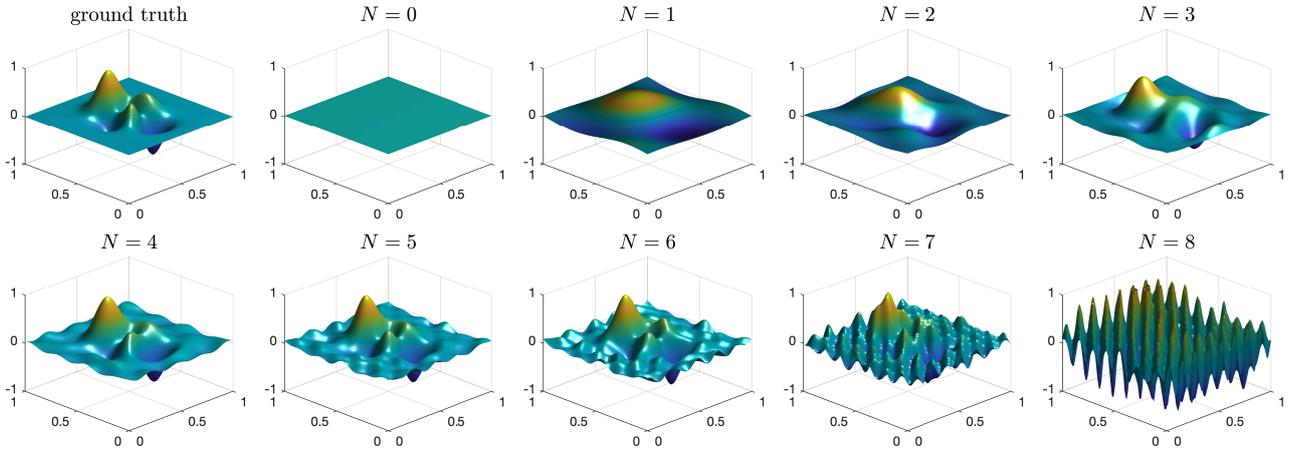

**Figure 8.** Surface plots of the reconstructions for the first surface profile with parameter values given in Table 2 (row 1) and different cut-off frequencies $N$.

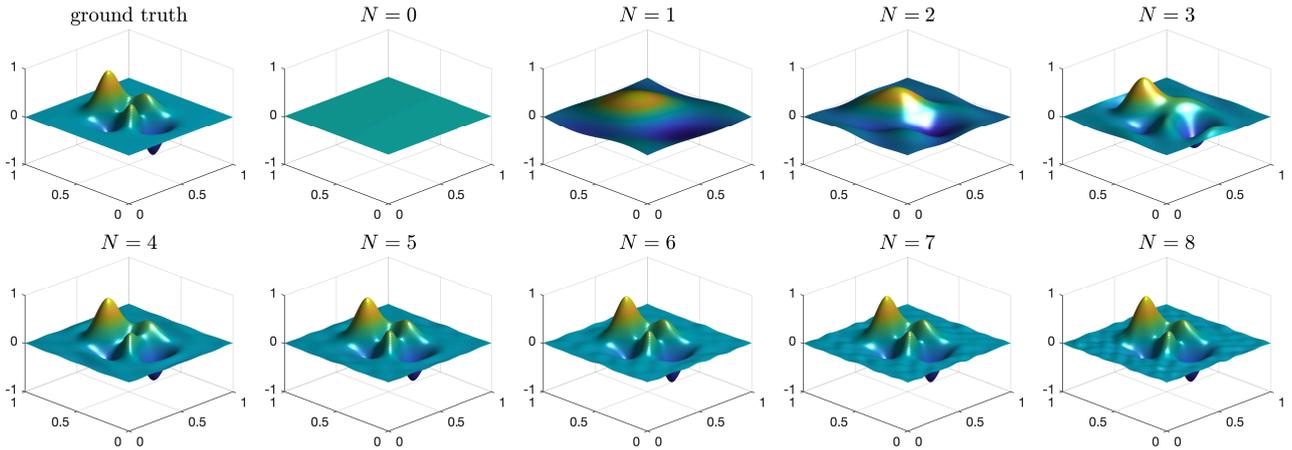

**Figure 9.** Surface plots of the reconstructions for the first surface profile with parameter values given in Table 2 (row 3) and different cut-off frequencies $N$.

According to (38), the reconstructions of all frequency modes are equally stable if $\rho = \kappa = -1$. Thus we expect a better stability for higher frequency modes as $\rho, \kappa$ approaches their ideal values. With $\rho = \kappa = -1 + 0.001i$ and roughly equal SNR, we obtain the results in Figure 9, which demonstrate better stability for higher frequency modes compared with the results with $\rho = \kappa = -1 + 0.01i$.

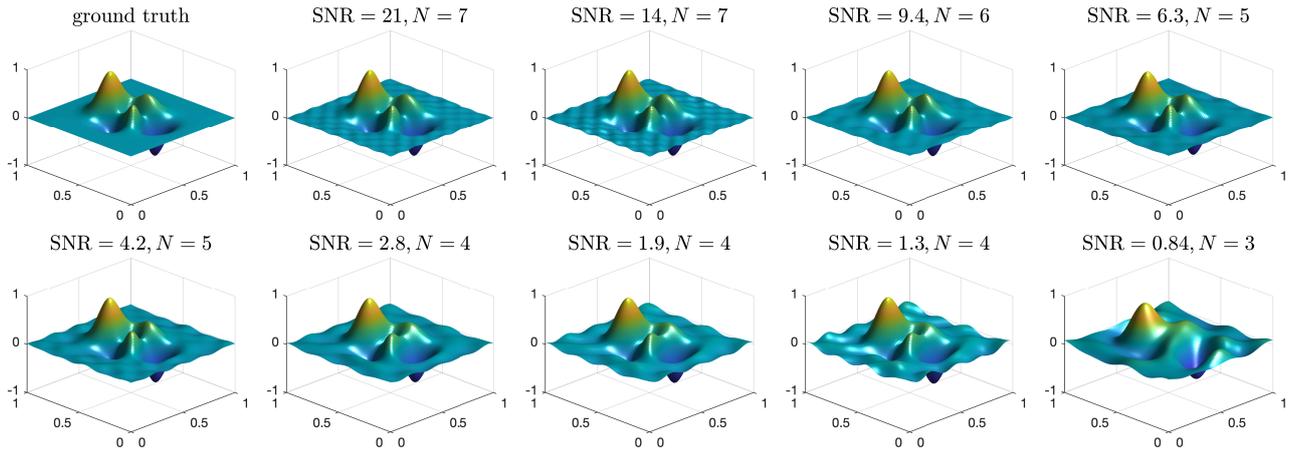

**Figure 10.** Pseudocolor plots of the reconstructions for the second surface profile using parameter values given in Table 2 (row 2), but with different values of SNR and the corresponding cut-off frequency $N$ determined from the discrepancy principle.

Figure 10 presents reconstructions conducted with the parameter values specified in Table 2 (row 2), but with different values of SNR and the corresponding cut-off frequency $N$ determined from the discrepancy principle. Again, the method exhibits robustness and over a wide SNR span, ranging from as low as 0.84 to as high as 21.

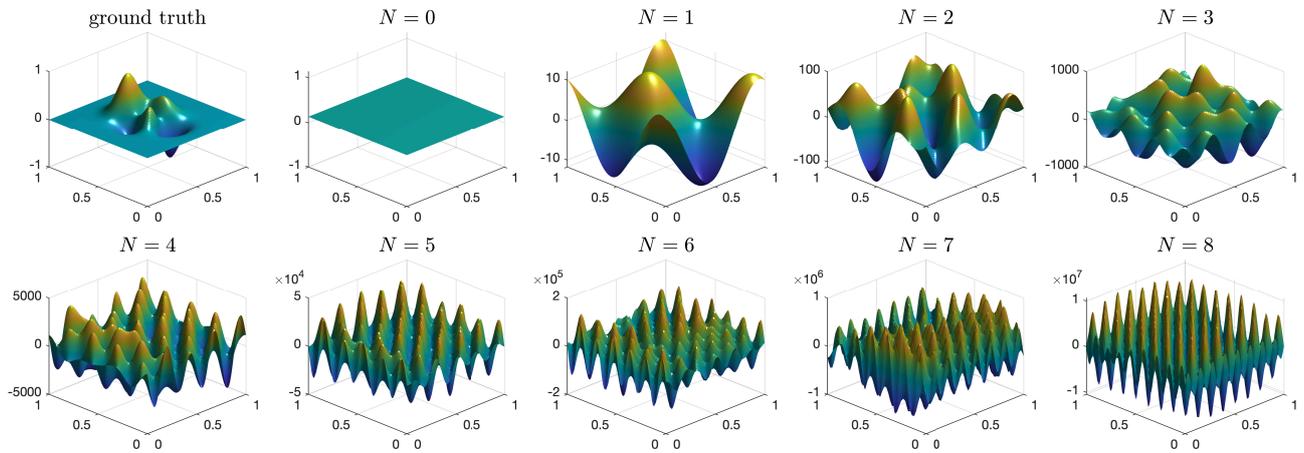

**Figure 11.** Surface plots of the reconstructions for the first surface profile with parameter values given in Table 2 (row 4) and different cut-off frequencies $N$.

Lastly, we show in Figure 11 the reconstructions when the slab is absent (by setting $\rho = \kappa = 1$) and with SNR = 9.3. Clearly, no frequency modes except $N = 0$ can be stably reconstructed.

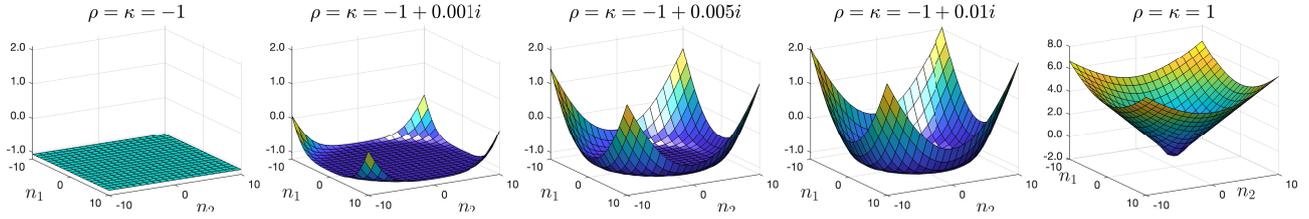

**Figure 12.** Plots of $\log_{10}|s_n|$ against $\boldsymbol{n}$ with different values of $\rho, \kappa$.

In Figure 12, we plot $\log_{10}|s_n|$ against $\boldsymbol{n}$ with different values of $\rho, \kappa$, where other parameters are given in Table 2 (row 2). With $\rho = \kappa = -1$, $|s_n|$ is a constant, indicating infinite resolution. If $\rho, \kappa$ deviates slightly from $-1$, then $|s_n|$ stays roughly constant and eventually increases exponentially as $\|\boldsymbol{n}\|$ increases, indicating increased instability for large frequencies. The farther away $\rho, \kappa$ deviates from $-1$, the earlier $|s_n|$ starts to increase. If $\rho = \kappa = 1$, then $|s_n|$ increases exponentially immediately from $\boldsymbol{n} = \boldsymbol{0}$. Note that $s_n$ is independent of $g$, thus the above observations apply to all surface profiles.

## 5.5 Experiment 3

Finally we consider the third surface profile as shown in Figure 2 (c).

|  | $\varepsilon$ | $a$ | $h$ | $\rho$ | $\kappa$ | $\sigma$ | SNR |
|---|---|---|---|---|---|---|---|
| Figure 13 | 0.001 | 0.1 | 0.1 | $-1 + 0.001i$ | $-1 + 0.001i$ | 0.0025 | 10.6 |
| Figure 14 | 0.01 | 0.1 | 0.1 | $-1 + 0.001i$ | $-1 + 0.001i$ | 0.031 | 10.6 |

**Table 3.** Parameter values used in Experiment 3.

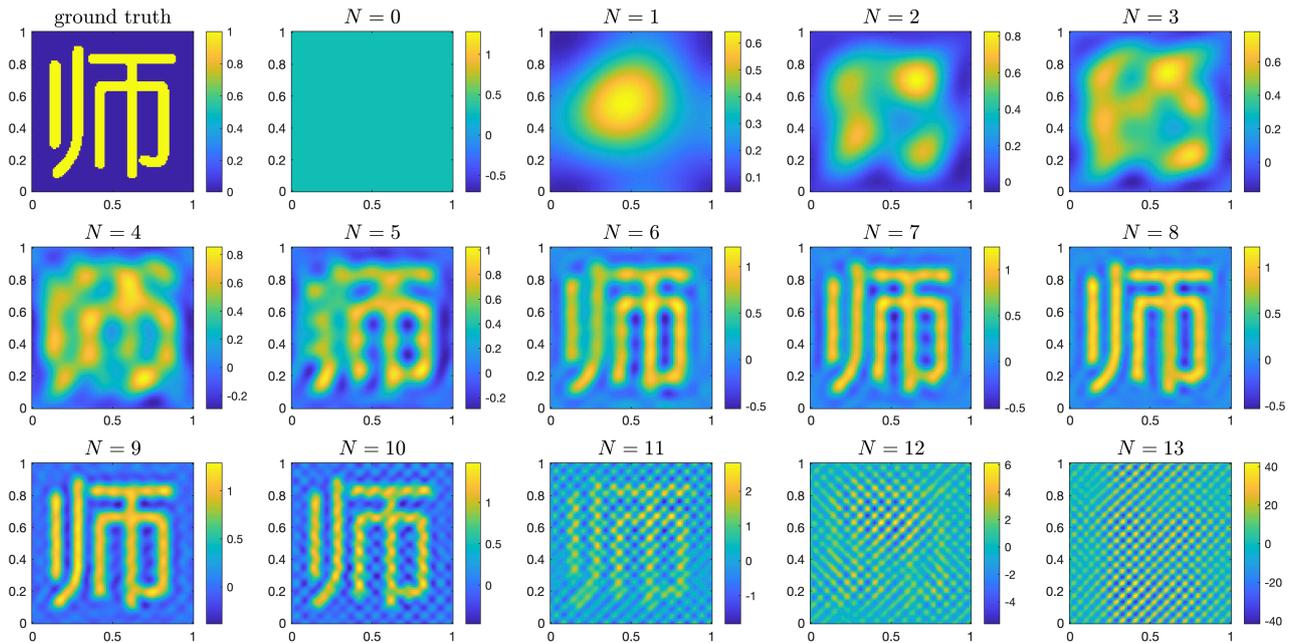

**Figure 13.** Surface plots of the reconstructions for the first surface profile with parameter values given in Table 3 (row 1) and different cut-off frequencies $N$.

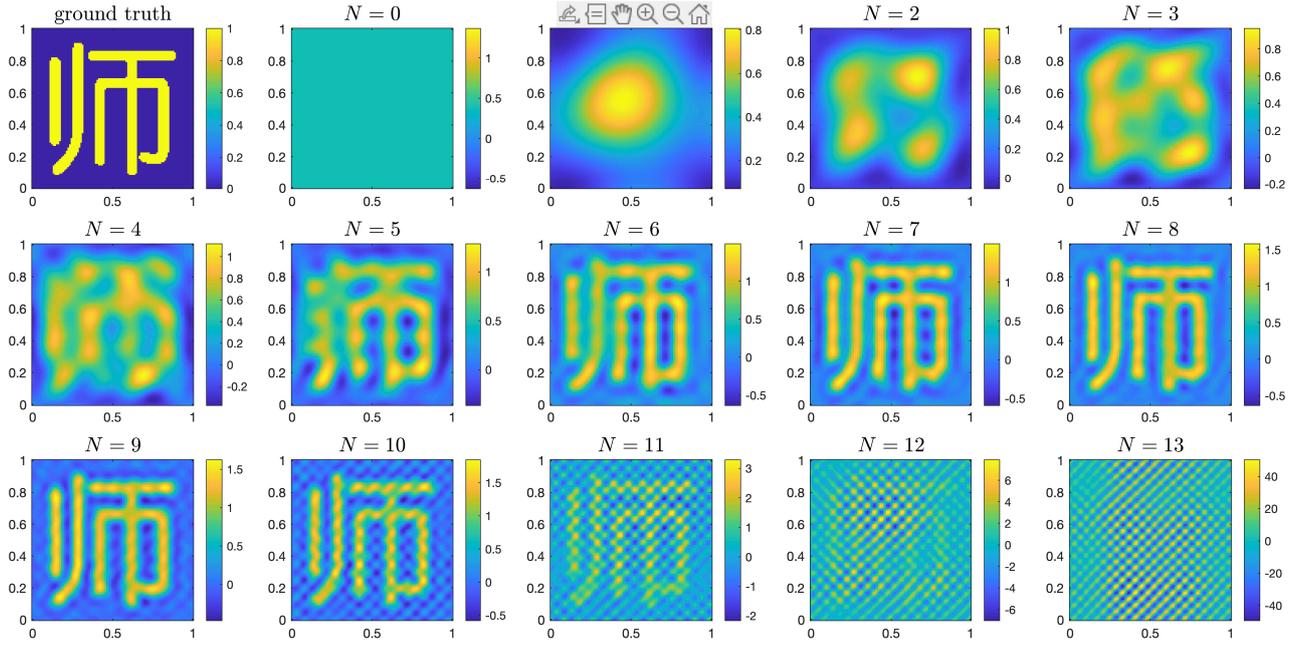

**Figure 14.** Surface plots of the reconstructions for the first surface profile with parameter values given in Table 3 (row 2) and different cut-off frequencies $N$.

With the parameter values given in Table 3 (row 1), we obtain the reconstructions in Figure 13. From visual inspection, the best results are approximately obtained with $N = 7, 8, 9$. Next, we increase the deformation parameter $\varepsilon$ from 0.001 to 0.01 while keeping the SNR approximately fixed. The results are shown in Figure 14, where the outline of the profile are as clear as those in Figure 13 at $N = 7, 8, 9$. However, the amplitude exhibit greater error as seen from the range in the colorbars, possibly due to the increasing of the linearization error as $\varepsilon$ increases.

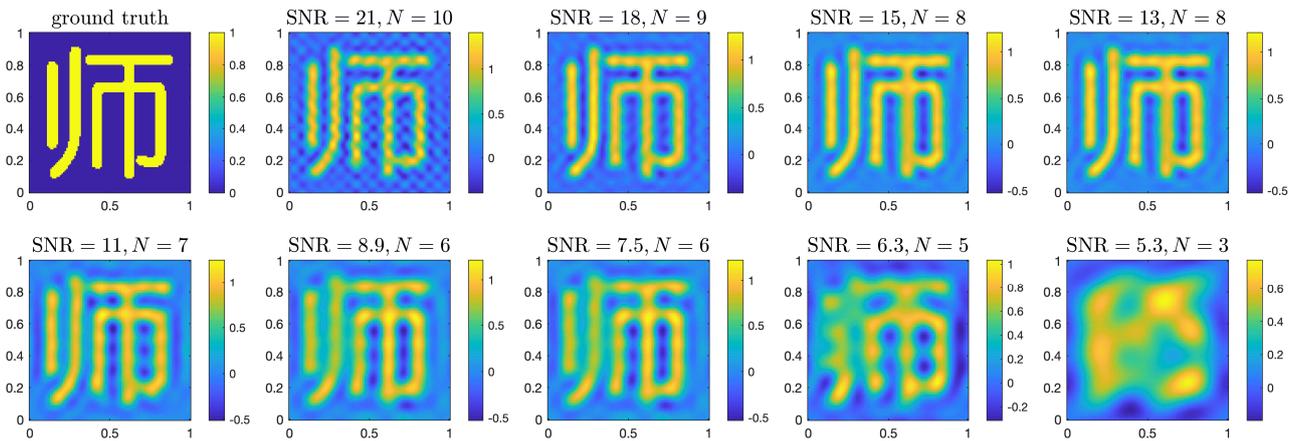

**Figure 15.** Pseudocolor plots of the reconstructions for the second surface profile using parameter values given in Table 3 (row 1), but with different values of SNR and the corresponding cut-off frequency $N$ determined from the discrepancy principle.

Figure 15 presents reconstructions conducted with the parameter values specified in Table 3 (row 1), but with different values of SNR and the corresponding cut-off frequency $N$ determined from the discrepancy principle with $c = 1.3$. Criterion to determine the optimal value of $c$ in the discrepancy principle is left for a future study.

# 6 Conclusion

We propose a model and numerical method for acoustic imaging of biperiodic surface with a slab of negative index material. It is formulated as an inverse scattering problem with measurement atop the slab from a single incident field. We first derive a boundary-interface value problem (BIVP) for the forward scattering problem and transform it to a rectangular domain through change of variables. Given a low-profile assumption on the surface, we take an asymptotic expansion of the wave field and derive a recursive system of BIVPs, which is solved to closed-form for the zeroth and first order terms. By dropping higher order terms in the expansion, we deduce a simple relation between the Fourier coefficients of the wave field and those of the surface profile, leading to a reconstruction formula for the inverse scattering problem. Preliminary analysis of the method shows an infinite resolution may be approximately achieved if the effective density and bulk modulus are both approximately $-1$ and the thickness of the prism is no less than the height of its lower surface. Numerical experiments for both smooth and nonsmooth surface profiles yield reconstructions with resolutions significantly beyond the diffraction limit, which can only be achieved by a much shorter measuring distance without the prism.

There are several open problems and extensions for the present study. The well-posedness of the forward problem is largely open and deserve a rigorous mathematical study. The uniqueness and stability of the inverse scattering problem is also worth study. Upon the establishment of those results, one may conduct convergence and error analysis of the proposed numerical scheme. If the low-profile assumption does not hold, then one may resolve to other numerical methods such as iterative or direct methods and still obtain significantly enhanced reconstructions. One can also take into consideration of the composite nature of the negative index metamaterial and conduct more realistic numerical or physical experiments. Finally one can extend the method to electromagnetic or elastic waves, and other geometrical settings such as inverse obstacle or medium scattering problems.

# Acknowledgments.

Our work was supported in part by the Guangdong Provincial Key Laboratory IRADS (2022B1212010006, R0400001-22).

# Appendix

In this appendix we study the DFT of white Gaussian noise. Denote by $g_\sigma$ the real-valued Gaussian distribution with mean zero and standard deviation $\sigma$. Suppose $u_i = v_i + iw_i$ is a complex-valued random variable such that both $v_i$ and $w_i$ follow the distribution $g_\sigma$ and mutually independent for $\boldsymbol{i} = (i_1, i_2)$, $0 \leqslant i_1 \leqslant I_1 - 1, 0 \leqslant i_2 \leqslant I_2 - 1$. Let $\boldsymbol{n} = (n_1, n_2)$, and

$$U_{\boldsymbol{n}} = \frac{1}{I_1 I_2} \sum_{i_1=0}^{I_1-1} \sum_{i_2=0}^{I_2-1} \eta_{\boldsymbol{ni}} u_{\boldsymbol{i}}$$

denote the DFT of $u_{\boldsymbol{i}}$, where

$$\eta_{\boldsymbol{ni}} = e^{-2\pi i \left(\frac{n_1 i_1}{I_1} + \frac{n_2 i_2}{I_2}\right)} = \lambda_{\boldsymbol{ni}} + i\mu_{\boldsymbol{ni}}.$$

Since $|\eta_{\boldsymbol{ni}}| = 1$, we see

$$\operatorname{Re} \eta_{\boldsymbol{ni}} u_{\boldsymbol{i}} = \lambda_{\boldsymbol{ni}} v_{\boldsymbol{i}} - \mu_{\boldsymbol{ni}} w_{\boldsymbol{i}}, \quad \operatorname{Im} \eta_{\boldsymbol{ni}} u_{\boldsymbol{i}} = \lambda_{\boldsymbol{ni}} w_{\boldsymbol{i}} + \mu_{\boldsymbol{ni}} v_{\boldsymbol{i}}$$

both follow the distribution $g_\sigma$. Moreover, one can verify

$$\begin{aligned}&\text{Cov}(\text{Re}\,\eta_{ni}u_i, \text{Im}\,\eta_{ni}u_i)\\ &= E[(\lambda_{ni}v_i - \mu_{ni}w_i)(\lambda_{ni}w_i + \mu_{ni}v_i)]\\ &= (\lambda_{ni}^2 - \mu_{ni}^2)E(v_i)E(w_i) + \lambda_{ni}\mu_{ni}[E(v_i^2) - E(w_i^2)] = 0\end{aligned}$$

Thus $\text{Re}\,\eta_{ni}u_i$ and $\text{Im}\,\eta_{ni}u_i$ are mutually independent for each $i$ and for fixed $n$. It follows that $\text{Re}\,U_n$ and $\text{Im}\,U_n$ are independent and both follow the distribution $g_{\sigma/\sqrt{I_1 I_2}}$.